\documentclass[%
 aip,
 amsmath,amssymb,
reprint, onecolumn %%%%% Comment onecolumn to override one column format %%%%%%%
]{revtex4-1}

\usepackage{graphicx}% Include figure files
\usepackage{dcolumn}% Align table columns on decimal point
\usepackage{bm}% bold math
\usepackage[utf8]{inputenc}
\usepackage[T1]{fontenc}
\usepackage{mathptmx}
\usepackage{etoolbox}

\usepackage{amsmath}
  \usepackage{paralist}
  \usepackage{amssymb}

\usepackage{color}

\def \beq{\begin{equation}}
\def \eeq{\end{equation}}

\def \ep{\epsilon}

\def \a{\alpha}

\makeatletter
\def\@email#1#2{%
 \endgroup
 \patchcmd{\titleblock@produce}
  {\frontmatter@RRAPformat}
  {\frontmatter@RRAPformat{\produce@RRAP{*#1\href{mailto:#2}{#2}}}\frontmatter@RRAPformat}
  {}{}
}%
\makeatother
\begin{document}

\preprint{AIP/123-QED}

\title[Eulerian droplet model]{ Initial boundary value problem for a system derived from Eulerian droplet model for air particle flow}
% Force line breaks with \\
\author{ Kayyunnapara Divya Joseph$^1$ }
 \email{divediv3@gmail.com}
\affiliation{ 
Department of Mathematics, \\
 Indian Institute of Science Education and Research Pune,\\
 Pune 411008, India%\\This line break forced with \textbackslash\textbackslash
}%
 \homepage{https://www.iiserpune.ac.in/research/department/mathematics/people/postdoctoral-fellows/post-doctoral-fellows/divya-joseph-kayyunnapara/525}

\date{\today}% It is always \today, today,

\begin{abstract}
 In this work, we study  the initial boundary value problem for a  non-strictly hyperbolic $2\times2$ system  of equations  in the quarter plane $x>0,t>0$  which is derived from Eulerian droplet model for air particle flow for velocity and volume fraction.  We show the existence of weak asymptotic solutions to the initial value problem  to the system  using  a regularisation,  by a vanishing viscosity method when the initial velocity is bounded measurable, the initial volume fraction is integrable and the boundary data are bounded measurable. Here we use a generalization of the Hopf-Cole transformation.  We also derive an explicit formula for the  weak solution  when the initial   data are functions of bounded variation,  the boundary datas are bounded and  locally in the class of Lipschitz continuous functions.  This construction involves  the Hopf-Lax formula for the boundary value problem for the Burgers equation and the product of a bounded variation function with derivative of another bounded variation function using  non-conservative Volpert product. \\
\end{abstract}
\maketitle
{ \flushleft
{\bf  AMS Subject Classification:} {35D30, 35F31, 35F50, 35L40}} \\
{\bf Keywords:}  non - strictly  hyperbolic system, vanishing viscosity, Hopf-Cole transformation,  Volpert product \\
{ \bf ORCID:} Kayyunnapara Divya Joseph https://orcid.org/0000-0002-4126-7882 \\
{ \bf E-Mail: } divediv3@gmail.com \\  
\footnote{ a) Department of Mathematics, Indian Institute of Science Education and Research Pune,  Dr. Homi Bhabha Road, Pashan, Pune 411008, Maharashtra, India.
}

\section{Introduction }
In this  work, we analyse a non - strictly hyperbolic system of conservation laws with linear damping given by, 
\begin{equation}
 \begin{aligned}
&\frac{\partial u}{ \partial t} + \frac{\partial }{ \partial x} \left( \frac{u^2}{2} \right) + \a(t) u=0, \\
&\frac{\partial v}{ \partial t} + \frac{\partial }{ \partial x} \left( u v \right) =0, 
\end{aligned}
\label{e1.1}
\end{equation}
in $x>0,t>0$, where $\a$ is a function of $t$. We study the initial boundary value problem in the quarter plane $x>0,t>0$, with initial data of the form 
\begin{equation}
 \begin{aligned}
u(x, 0)= u_0(x), \,\,\, v(x, 0)= v_0(x),\,\,\, x>0,
\end{aligned}
\label{e1.2}
\end{equation}
various boundary  conditions for $u$ at $x=0$ and  the mass condition for $v$
\begin{equation}
u(0,t)=u_B, \,\,\int_0^\infty v(y,t)dy =v_B(t), t>0.
\label{e1.3}
\end{equation}
We also study \eqref{e1.1} in $x>0,t>0$ with Dirichlet boundary condition 
\begin{equation}
u(0,t)=u_B(t), \,\,v(0,t) =v_B(t), t>0,
\label{e1.4}
\end{equation}at $x=0$.
 We see that the boundary condition \eqref{e1.3} or \eqref{e1.4} cannot be prescribed in the strong sense.

When $\a>0$, a constant, the system \eqref{e1.1} is derived from the Eulerian droplet model  for air particle flow, see  \cite{Keita}. Here $\a$ is the drag coefficient between the air and the particles and we consider the case when the drag  coefficient depends on time. Also the components  $v \geq 0$ and $u$ are the volume fraction and the velocity of the particles respectively. The system \eqref{e1.1} has repeated eigenvalue $u$ of multiplicity $2$ with one dimensional eigenspace, so it is a  non-strictly hyperbolic system.

When $\a =0$, the  first equation in system \eqref{e1.1} is commonly called the Burger's equation, that  was first derived in 1915 by Bateman  \cite{Bateman}, later reinvented by Burgers  \cite{Burgers} in 1948.  It is useful in understanding various phenomenon like shock waves in gas dynamics, see the work of Hopf \cite{Hopf}. In fact
when $\a =0,$ $u$ is constant along the characteristics defined by $\frac{dx}{dt}=u$  and the characteristics are straight lines.  For non-increasing initial data, the characteristics starting at different points of the initial line meet  at some later time and shocks are generated. 

The system \eqref{e1.1}  with $\alpha=0$ is used in modelling the evolution of density inhomogeneities in matter in the universe, see the work of  Shandarin and Zeldovich \cite{Shandarin}. There are more recent applications involving kinetic models of stochastic production flows such as the flows of products through a factory or supply chain. These kinetic models when expanded into deterministic moment equations lead to the system \eqref{e1.1} with $\a=0$, see \cite{a1, a2,  Forestier}. The particles  moving along the characteristics lines collide,  stick together and become massive particles. Due to this concentration phenomenon, $v$ generally becomes  a measure. 
Systems of this type  were analyzed by many authors  see the work of S. N. Gurbatov and A.I. Saichev \cite{Saichev}. S. Albeverio and  V. M. Shelkovich \cite{Albeverio} worked on general system including the zero pressure gas dynamics by finding the specific  type of solutions for the Cauchy problem with piecewise-smooth
initial data. They proved that these are $\delta$ waves solutions,  which are related to
the concentration process on the surface which carries the singularities.  In  \cite{LeFloch, Tan} the well posedness of the Cauchy problem to  the system \eqref{e1.1} with $\alpha=0$ was  analyzed. When $\a\neq 0$ we have
the same phenomenon, but the characteristics are not straight lines and  $u$ is not constant along the characteristics. Indeed $\frac{du}{dt}$ along the characteristics is $- \a u$.  

 The pure initial value problem for \eqref{e1.1} was studied by many authors. Richard \cite{Richard} studied the case of the Riemann initial value problem  and Divya \cite{Divya} studied the case for general initial data.  In these works, a viscosity approximation  with a variable viscosity coefficient and a small parameter $\epsilon>0$ was used. By a generalized Hopf-Cole transformation,  they reduced the system into an initial value problem for the heat equation and  got an explicit formula for the  approximate solution. By letting $\ep \rightarrow 0$,  Richard \cite{Richard} got an  explicit formula for the Riemann problem to \eqref{e1.1}.  Divya \cite{Divya} used the weak asymptotic solution approach of Maslov adapted by Danilov, Omelyanov and  Shelkovich \cite{Danilov}, Albeverio and Shelkovich \cite{Albeverio} for their study of the initial value problem for a class of problems where measure solutions appear.
 Divya \cite{Divya}  constructed  the explicit weak asymptotic solution for general initial data. Further, when the initial data  $(u_0,v_0)$  is in the space of functions of bounded variation, by taking the limit as $\epsilon \rightarrow 0$, an explicit solution was constructed  in the sense of measures. 
 
In our present work,  we study the initial boundary value problem to system \eqref{e1.1} with different types of boundary conditions.  Here, we prove existence of  a weak asymptotic solution  for general initial and boundary data. We are not able to find the solution explicitly. Even in the case  when the explicit formula is available, passing to the limit is difficult due to the formation of boundary layers, as the limit may not satisfy the strong form of boundary condition, see \cite{Joseph3}, the analysis for the Burgers equation. So, we use another approach using a transformation to reduce the case  to that when $\alpha=0$ and we construct explicitly the measure solution to the initial value problem with a suitable weak form of boundary conditions. The crucial step is to use the formula for the initial boundary value problem derived in \cite{Joseph3,Joseph5} for the Burgers equation. For this construction, we need to assume  that the initial and boundary data are regular. In particular, we assume that $(u_B(t),v_B(t))$ is locally Lipschitz continuous on $[0, \infty)$. For special cases,  the solution obtained by passing to the limit in the weak asymptotic solution and the one obtained by the transformation agree.

 Rest of the paper  is organised as follows.  Section II contains the construction of a  weak asymptotic solution for the boundary value problem to  \eqref{e1.1}-\eqref{e1.3}. In  Section III, we construct an explicit solution to the initial  boundary value problem to  \eqref{e1.1}-\eqref{e1.3} using a transformation. We consider the cases when $\a =0$ in  Section A and when $\a(t) \neq 0$  in  Section B separately. We also   construct an explicit solution to the initial  value problem  for \eqref{e1.1}-\eqref{e1.3} in  Section C. The last section contains the conclusions of our study.

\section{Construction of weak asymptotic solution}

In this section we construct the weak asymptotic solution for the boundary value problem \eqref{e1.1}-\eqref{e1.3}.  First we explain the notion of weak asymptotic solution following 
the work of Albeverio and Shelkovich \cite{Albeverio} which was for the initial value problem. In our case we need to incorporate boundary conditions also in the weak formulation of the boundary value problem.
{ \flushleft
 {\bf Definition 2.1:}}\\
 A family of smooth functions $\left( u^{\epsilon}, v^{\epsilon} \right)_{\epsilon >0}$ is called a {\bf weak asymptotic solution} to the system
\eqref{e1.1} with initial conditions \eqref{e1.2} and boundary conditions \eqref{e1.3} provided: \\
\begin{itemize}
\item i) As $\epsilon \rightarrow 0$ we have, for all  $\phi \in C_c^{\infty}(0, \infty),$ 
\begin{equation}
\begin{aligned}
\int_0^{\infty}{ ( u^\epsilon_t(x, t) + u^\epsilon(x, t)  u^\epsilon_x(x, t) +   \a(t) u} ) \phi(x) dx &=o(1),\\
\int_0^{\infty}{ ( \rho^\epsilon_t(x, t)  + {(v^\epsilon(x, t)  u^\epsilon(x, t) )}_x ) \phi(x)  \,\,\, dx} &= o(1),\\  
\end{aligned}
\label{e1.5}
\end{equation}
 The estimates \eqref{e1.5} are required to hold uniformly on $[0, T]$ for each $T$.\\
\item ii) As $\epsilon \rightarrow 0,$ for all  $\phi \in C_c^{\infty}(0, \infty), $ we have, 
\begin{equation}
\begin{aligned}
\int_0^{\infty}{ ( u^\epsilon(x,0)-u_0(x) ) \phi(x)  dx} &=o(1), \int_0^{\infty} {( v^\epsilon(x,0)-v_0(x) ) \phi(x)  dx} =o(1),\\
\int_0^{\infty} {( u^\epsilon(0,t)-u_B(t)) \phi(t)  dt} &=o(1), \int_0^{\infty}{ ( \int_0^\infty v^\epsilon(y,t) dy-v_B(t)) \phi(t)  dt } =o(1).
\end{aligned}
\label{e1.6}
\end{equation}
\end{itemize}

{ \flushleft
 The boundary conditions are not satisfied strongly in the passage to the limit. In fact for special cases where we can actually compute the limit $\lim_{\epsilon \rightarrow 0}(u^\epsilon,v^\epsilon)=(u,v)$, we will see that a weak form of boundary conditions are satisfied namely }
\begin{equation}
u(0,t) \in E(u_B(t)),\,\,v(0,t)= v_B(t), \,\, if \,\,\, u(0,t)>0,\,\, \text{ a.e.} \,\,\, t>0
\label{e1.7}
\end{equation}
 For $u_B$ the admissible set $E(u_B)$ is defined by
\begin{equation}
E(u_B) =  \begin{cases} \displaystyle
      {(-\infty,0]) ,\,\,\,\text{if} \,\,\,u_B \leq 0}
\\\displaystyle
        {\{u_B\}\cup (-\infty, -u_B],\,\,\,\text{if} \,\,\, u_B>0}.
\end{cases}
\label{e1.8}
\end{equation}
We prove the following theorem.\\
{ \flushleft
{\bf  Theorem 2.2: }} Let $u_0$ , $u_B$  and $v_B$ be bounded measurable functions on $[0,\infty)$ and $v_0$ be integrable function on $[0,\infty)$. Then there is a weak asymptotic solution of \eqref{e1.1} with initial condition \eqref{e1.2} and boundary condition \eqref{e1.3}.\\

{ \flushleft
The main ingredient of the proof of the theorem is a regularisation  using a variable vanishing viscosity approximation, introduced by Richard \cite{Richard} in his study of  the Riemann problem. }
\begin{equation}
 \begin{aligned}
&\frac{\partial u}{ \partial t} + \frac{\partial }{ \partial x}  \left( \frac{u^2}{2} \right) + \a(t) u =  {\ep} \,\, \left(e^{-\int_0^t \alpha(s') ds'} \right) \,\, \frac{{\partial}^2 u}{ \partial x^2}, \\
&\frac{\partial v}{ \partial t} + \frac{\partial }{ \partial x} \left( u v \right) =  {\ep} \,\, \left(e^{-\int_0^t \alpha(s') ds' } \right)\,\,\frac{{\partial}^2 v}{ \partial x^2},
\end{aligned}
\label{e2.1}
\end{equation}
in the space - time domain $\{(x,t) : x>0, t>0\}$ with initial conditions
\begin{equation}
 \begin{aligned}
u(x, 0)= u_0(x), \,\,\, v(x, 0)= v_0(x), x >0
\end{aligned}
\label{e2.2}
\end{equation}
and boundary conditions 
\begin{equation}
u(0,t)=u_B(t),\,\,\,\int_0^\infty v(y,t) dy=v_B(t).
\label{e2.3}
\end{equation}
We  construct an explicit solution to this Cauchy problem. \\

{ \flushleft 
{ \bf Theorem 2.3: }} Let $u_0 \in L^\infty( 0,\infty)$ and $v_0$ be a function of bounded variation which is in $L^1( 0,\infty)$.  
Then there exists a unique solution  $(u^{\epsilon},v^{\epsilon})$  to \eqref{e2.1} - \eqref{e2.3}  which is $C^\infty$ in $x>0,t>0$, $v^{\epsilon}$ with $v^{\epsilon} = V_x$  an  integrable function of $x$ for each $t>0$. 
$u^{\epsilon}$ and $V=-\int_x^\infty v(y,t) dy$ are bounded uniformly in $\epsilon$.
{ \flushleft
{ \bf Proof: }}  To prove the theorem, we introduce
\begin{equation}
A(s)=e^{\int_0^t \alpha(s) ds} ,\,\,\tau=\tau(t)=\int_0^t \frac{ds}{A(s)}.
\label{e2.4}
\end{equation}
 Clearly $\tau : [0,\infty) \rightarrow  \left[0,\int_0^\infty \frac{ds}{A(s)} \right)$ is one-one and so an invertible  map. We use the notation $\tau(t)$ and $t(\tau)$  in \eqref{e2.4}  as inverses of each other.\\
 
  First we claim that, if $(p(x,\tau),q(x,\tau))$ solves
\begin{equation}
 \begin{aligned}
 p_\tau +p p_x =\epsilon p_{xx},\,\,q_\tau +(p q)_x = q_{xx}
\end{aligned}
\label{e2.5}
\end{equation}
in $x>0,\tau>0$ with initial data 
\begin{equation}
p(x, 0)=u_0(x) \,\,\,q(x, 0)=v_0(x) 
\label{e2.6}
\end{equation}
and boundary condition 
\begin{equation}
p(0,\tau)=A(t(\tau)) u_B(t(\tau)),\,\,\int_0^\infty q(y,\tau)dy =v_B(t(\tau))
\label{e2.7}
\end{equation}
then 

\begin{equation}
 \begin{aligned}
(u^\epsilon,v^\epsilon)=\left(\frac{1}{A(t)} p(x,\tau(t)),q(x,\tau(t))\right)
\end{aligned}
\label{e2.8}
\end{equation} 
solves \eqref{e2.1}  in $x>0,t>0$ with initial condition \eqref{e2.2} and boundary condition \eqref{e2.3}. The claim on the initial and boundary conditions are obvious.
Now coming to the partial differential equations we  compute derivatives of $A(t)u^\epsilon(x,t)=p(x,\tau)$ and $v^\epsilon(x,t)=q(x,\tau)$ as follows.
\[
\begin{aligned}
&A'(t) u^\epsilon(x,t) +A(t) u^\epsilon_t (x,t)=p_\tau \frac{d \tau}{dt}, \,\,\, A(t) u^\epsilon_x(x,t)=p_x(x,\tau),  \,\,\,  A(t) u^\epsilon_{xx}(x,t) =p_{xx}(x,\tau),\\
&v_t(x,t)=q_\tau (x,\tau)\frac{d \tau}{dt},  \,\,\,   v_x(x,t)=q_x(x,\tau),  \,\,\,   v_{xx}=q_{xx}(x,\tau).
\end{aligned}
\]
Using $A'(t)=\alpha(t) A(t)$ and  $\frac{d \tau}{dt}=\frac{1}{A(t)}$, in these equations we get
\[
\begin{aligned}
&p_\tau +p p_x -\epsilon p_{xx}=A^2(t)  \left[u_t+\alpha(t) u +u u_x -\epsilon \frac{1}{A(t)} u_{xx}\right]\\
&q_\tau + (p q)_x -\epsilon q_{xx} =A(t) \left[v_t +(uv)_x -\frac{1}{A(t)} v_{xx}\right]
\end{aligned}
\]
Since $A(t)\neq 0$  our claim follows.

Next, note that  if $(P,Q)$ is solution to 
\begin{equation}
 \begin{aligned}
 P_\tau +\frac{P_x P_x}{2} =\epsilon P_{xx},\,\,Q_\tau +(P_xQ_x) = \epsilon Q_{xx}
\end{aligned}
\label{e2.9}
\end{equation}
with initial conditions
\begin{equation}
 \begin{aligned}
 P(x,0)=\int_0^x u_0(s) ds,\,\,\,Q(x, 0)=-\int_x ^\infty v_0(s) ds\
 \end{aligned}
\label{e2.10}
\end{equation}
and boundary condition 
\begin{equation}
P_x(0,\tau)=A(t(\tau)) u_B(t(\tau)),\,\,Q(0,\tau) =-v_B(t(\tau))
\label{e2.11}
\end{equation}
 then
\begin{equation}
(p,q) = (P_x,Q_x)
\label{e2.12}
\end{equation}
solves \eqref{e2.5} -\eqref{e2.7}.\\

We use the Hopf-Cole transformation \cite{Hopf} and \cite{Joseph4}
\begin{equation}
 \begin{aligned}
P= - 2 {\ep} \log{S^{\ep}}, \,\,\,
Q = \frac{C^{\ep}}{S^\ep}. 
\end{aligned}
\label{e2.13}
\end{equation} 
to reduce the problem \eqref{e2.9}-\eqref{e2.11} to  the system,  
\begin{equation}
 \begin{aligned}
C^{\ep}_{\tau}  =  {\ep}  C^{\ep}_{ xx}, \,\,\,
S^{\ep}_{\tau}  =  {\ep}  S^{\ep}_{ xx}  
\end{aligned}
\label{e2.14}
\end{equation}
with initial conditions
\begin{equation}
S^\ep(x,0)=\exp{ \left\{{-\frac{\int_0^x u_0(s)ds}{2 \ep}} \right\}},\,\,\,C^\ep(x,0)=-\int_x^\infty v_0(s) ds \,\,  \exp{\left\{{-\frac{\int_0^x u_0(s)ds}{2 \ep}} \right\}}
\label{e2.15}
\end{equation}
and boundary condition

\begin{equation}
2\epsilon S_x^\ep(0,\tau)+A(t(\tau))u_b(t(\tau)) S(0,\tau)=0,\,\,  C(0,\tau)=-v_B(t(\tau)) S(0,\tau)
\label{e2.16}
\end{equation}

By linear theory, there exists a unique solution $({S^\ep}, {C^{\ep}})$, to \eqref{e2.14} - \eqref{e2.16}. By the regularity theory of the heat equation, they are $C^\infty$ in $x>0,\tau>0$.
 Now using \eqref{e2.8}, \eqref{e2.12} and \eqref{e2.13},  we can express $(u^\ep, v^\ep)$ in terms of  ${C}^\ep, {S}^\ep$ and their first
$x$ derivatives. So $(u^\epsilon, v^\epsilon)$ are $C^\infty((0,\infty) \times (0,\infty))$.

 Now we get the $L^\infty$ estimates for $u^\epsilon$ and $V^\epsilon$. Applying the maximum principle to this differential equation for $p$,  we get
\[
||p||_{L^\infty((0,\infty) \times[0,\tau_0))}\leq ||u_0||_{L^\infty(0,\infty)} +||u_B(t)||_{L^\infty(0,t(\tau_0))},
\]
for any $\tau_0>0$. So, for $T>0$ we get,
\begin{equation}
||{u}^\epsilon||_{L^\infty(( 0,\infty) \times [0,T))} \leq e^{||\alpha || _{L^1[0,T]}} \left( ||u_0||_{L^\infty(0,\infty)} +||u_B(t)||_{L^\infty(0,T]}\right)
\label{e2.17}
\end{equation}

Also  $Q^\epsilon$ satisfies the equation
\[
\begin{aligned}
&Q_t\tau+p Q_x =\epsilon Q_{xx}, x>0,\,\,\tau>0\\
&Q(x,0)=-\int_x^\infty v_0(y) dy, \,\,x.>0,\,\,Q(0,t)=-v_B(t)
\end{aligned}
\]
Again by the maximum principle we have for any $\tau_0>0$,
\begin{equation}
||{Q}^\epsilon||_{L^\infty ((0,\infty) \times [0,\tau_0))}\leq ||v_0||_{L^1(0,\infty)} +||v_B||_{L^\infty(0,t(\tau_0))} 
\label{e2.18} 
\end{equation}
Now note that $V(x,t)=-\int_x^\infty v(y,t) dy =Q(x,\tau(t))$, as they satisfy the same equations and boundary conditions. So we have the estimate
\[
||V^\epsilon||_{L^\infty ((0,\infty) \times [0,T))}\leq ||v_0||_{L^1(0,\infty)} +||v_B||_{L^\infty(0,T)}
\]
for any $T>0$.
Uniqueness of solution $u^\epsilon$ is standard as it is a scalar parabolic equation with smooth initial data. The second equation is a scalar linear parabolic equation for $v^\epsilon$ with 
a smooth coefficient, with initial and boundary conditions whose solution is also unique by the classical theory of parabolic equations.
{ \flushleft
{\bf Proof of theorem 2.2:} }\\
 Now we consider the general initial and boundary data to show the existence of a weak asymptotic solution.
 
  First, we regularise the initial and boundary data using a cut off  function near $0$ and then regularise it by a convolution in the scale $\epsilon$ as follows.
Let $\chi_{[2\epsilon,\infty)}$ be the characteristic function on $[ 2 \epsilon, \infty)$ and $\eta_\epsilon$ be the usual Friedrichs mollifier in one space dimension. 
Let
\begin{equation}
\begin{aligned}
&u^\epsilon_0(x)=(u_0 \chi_{[2\epsilon,\infty)}*\eta^\epsilon)(x),\,\,\, v^\epsilon_0(x)=(v_0 \chi_{[2\epsilon,\infty)}*\eta^\epsilon)(x),\\
&u^\epsilon_B(t)=(u_B \chi_{[2\epsilon,\infty)}*\eta^\epsilon)(t),\,\,\, v^\epsilon_B(t)=(v_B \chi_{[2\epsilon,\infty)}*\eta^\epsilon)(t).
\end{aligned}
\label{2.19}
\end{equation}
Now let $(\tilde{u}^\epsilon,\tilde{v}^\epsilon)$ be the solution to \eqref{e2.1} with initial conditions
\[
\tilde{u}^\epsilon(x,0)=u^\epsilon_0(x),\,\,\,\tilde{v}^\epsilon(x,0)=v^\epsilon_0(x), x>0
\]
and boundary conditions 
\[
\tilde{u}^\epsilon(0,t)=u^\epsilon_B(t),\,\,\,\int_0^\infty \tilde{v}^\epsilon(y,t)dy=v^\epsilon_B(t), t>0
\]
 constructed by the Theorem 2.3. Then clearly $(\tilde{u}^\epsilon,\tilde{v}^\epsilon)$  is a $C^\infty$ function. Also since 
\[
\begin{aligned}
&||u^\epsilon_0||_{L^\infty(0,\infty)} \leq ||u_0||_{L^\infty(0,\infty)},\,\,||u^\epsilon_B||_{L^\infty(0,\infty)} \leq ||u_B||_{L^\infty(0,\infty)}\\
&||v^\epsilon_0||_{L^1(0,\infty)} \leq ||v_0||_{L^1(0,\infty)},\,\,||v^\epsilon_B||_{L^\infty(0,\infty)} \leq ||v_B||_{L^\infty(0,\infty)},
\end{aligned}
\] 
using the estimates \eqref{e2.17} and \eqref{e2.18}, we get a constant $C(T)$ independent of $\epsilon$ such that
 \[
 ||\tilde{u}^\epsilon||_{L^\infty((0,\infty) \times [0,T])} \leq C(T)\,\,\,,
||\int_x^\infty \tilde{v}(y,t)||_{L^\infty((0,\infty) \times [0,T])} \leq C(T).
\]
Now consider any test function $\psi \in C_c^\infty(0,\infty)$ and the following integrals,
\begin{equation}
 \begin{aligned}
\int_0^\infty( {\tilde{u}^\epsilon}_t +  {\tilde{u}^\epsilon} {\tilde{u}^\ep}_x+ \a(t) \tilde{u}^\epsilon) \psi dx &= \int_0^\infty {\ep} e^{- \int_0^t \a(s)ds}  {\tilde{u}^\epsilon}_{xx} \psi(x) dx\\
&=\int_0^\infty{\ep} e^{- {\int_0^t \a(s) ds}} {\psi}_{xx} \tilde{u}^\epsilon(x,t) dx, \\
\int_0^\infty( {\tilde{v}^\epsilon}_t + (\tilde{u}^\epsilon \tilde{v}^\epsilon)_x)\psi(x) dx &=  \int_0^\infty{\ep} e^{- \int_0^t\a(s) ds}  {\tilde{v}^\epsilon}_{xx} \psi(x)dx \\
&=\int_0^\infty{\ep} e^{- \int_0^t \alpha (s)ds}  { \psi(x)}_{ xxx} \left(\int_x^\infty \tilde{v}^\epsilon(y,t) dy \right) dx.
\end{aligned}
\label{e2.20}
\end{equation}
Clearly the  terms on the  extreme right of \eqref{e2.20} goes to zero for $\epsilon \rightarrow 0$ unifomly in $t \in [0,T]$, for every $T>0$. Also from \eqref{2.19}, by the property of the convolution we have,
\[
\tilde{u}^\epsilon(x,0)-u_0 (x)\rightarrow 0,\,\,\tilde{v}^\epsilon(x,0) -v_0(x)\rightarrow 0, \,\,\tilde{u}^\ep(0,t) -u_B (t)\rightarrow 0,\,\,\int_0^\infty \tilde{v}^\ep(y,t)dy-v_B(t) \rightarrow 0, 
\]
in distribution as $\epsilon \rightarrow 0,$ which is \eqref{e1.6}. This proves that $(\bar{u}^\epsilon,\bar{v}^\epsilon)$ is a weak asymptotic solution of \eqref{e1.1} - \eqref{e1.3}.

\subsection{\label{sec:level2}Dirichlet boundary condition}
 In this section, we work on the analysis for the Dirichlet boundary value problem  to \eqref{e1.1} in $x>0,t>0$, where $\a$ is a function of $t$. We study the initial boundary value problem for \eqref{e1.1} in the quarter plane $x>0, t>0$, with initial data of the form  \eqref{e1.2} and  Dirichlet boundary condition \eqref{e1.4}.  The analysis is similar. We use the same parabolic regularisation \eqref{e2.1} with regularised initial and boundary data  in \eqref{e1.2}, \eqref{e1.4}  as before. We use the same transformations as in equations \eqref{e2.8},\eqref{e2.12} and \eqref{e2.13} to get the equations \eqref{e2.14} for $C^\epsilon$ and $S^\epsilon$. The initial conditions for $C^\epsilon$ and $S^\epsilon$ and boundary condition for $S^\epsilon$ is also the same as before. The only difference is that the boundary condition  $v(0,t)=v_B(t)$ in
\eqref{e1.4}  is transformed for the linearised problem for $(S,C)$ in the following way
\begin{equation}
2\ep C_x^\ep(0,\tau)+A(t(\tau)) \,\,\, u_B(t(\tau)) \,\,\, C^\ep(0,\tau)=2 \epsilon \,\,\, v_b(\tau(t)) S^\ep(0,\tau).
\label{e2.24}
\end{equation}
 To get this boundary condition for $C^\ep$, we just observe that
\[
q=Q_x=(\frac{C}{S})_x=\frac{C_x}{S}-\frac{C}{S} \frac{S_x}{S}=\frac{C_x}{S}-\frac{C}{S} \frac{-p}{2 \ep}=\frac{2 \ep C_x + p C}{2 \ep S}.
\]
This gives
\[
2 \ep C_x + p C = 2 \ep q S
\]
and the boundary condition   \eqref{e2.24} for $C^\epsilon$ follows as 
\[
q(0,\tau)=v_B(t(\tau)), p(0,\tau)=A(t(\tau)) u_B(t(\tau)).
\]

\section{Solution by the method of transformation - the case when initial and boundary data are regular}

 In this section, we  derive an explicit formula  for the solution of the  system \eqref{e1.1} with initial condition  \eqref{e1.2} and boundary conditions  \eqref{e1.3}.  The analysis in this case is much more difficult than the case for the pure initial value problem.  For the initial value problem in $\{(x,t) : x\in R, t>0 \}$,  an explicit formula for the weak asymptotic solution  was derived in \cite{Divya} for the general bounded measurable initial  data $u_0$, integrable $v_0$, using the regularisation \eqref{e2.1}, when $\alpha$ is a constant. Further,  when  initial data in the space of bounded variation,  an explicit formula was derived for the vanishing viscosity limit.  For the initial boundary value problem  in $\{(x,t) : x>0, t>0 \}$, we do not have an explicit formula for the viscosity approximation \eqref{e2.1}. Even for special cases where explicit formula is available,  boundary layers appear as the limit does not satisfy the given boundary condition  for the viscosity approximation, see  \cite {Divya1} and references therein. So here, we use another approach
to solve the initial boundary value problem  explicitly.  We introduce a new set of  independent variable $(x,\tau)$ and dependent variable $( p(x,\tau,q(x,\tau))$ defined by  \eqref {e2.4} and

\begin{equation}
 p(x,\tau)=A(t(\tau)) u(x,t(\tau)), \,\,\, q(x,\tau)=v(x,t(\tau))
\label{3.3}
\end{equation}

The motivation for this transformation comes from the work  \cite{Joseph1}, where a single equation with a general damping term was studied.  For the weak asymptotic solution, the initial and boundary data are regularised by convolution. But for the explicit formula derived here,  we need to assume that  $u_0$  is a function of  bounded  variation, $v_0$  an
integrable function which is of  bounded variation in $x >0$ and  $u_B(t)$, $v_B(t)$  are   locally 
Lipschitz continuous functions   in   $t>0$.
{ \flushleft
{\bf Proposition 3.1: }} $(p,q)$ is a a smooth solution to 
\begin{equation}
\begin{aligned}
&p_\tau + \left(\frac{p^2}{2} \right)_x=0, \\
&q_\tau + \left(p q \right)_x=0
\end{aligned}
\label{e3.1}
\end{equation}
 iff  $(u,v)=(A(t)^{-1} p(x,\tau(t)),q(x,\tau(t)))$ is a smooth solution to  \eqref{e1.1}.\\
{\bf Proof :} An easy computation shows that under \eqref{3.3},  $\frac{d \tau}{dt}=\frac{1}{A(t)}$ and so
\[
\begin{aligned}
&p_\tau = A'(t) \frac{dt}{d \tau} u+ A(t) \frac{dt}{d \tau} u_t=A^2(t) (u_t +\alpha(t) u),\\
&q_\tau =v_t \frac{dt}{d \tau} = A(t) v_t\\
&\left(\frac{p^2}{2}\right)_x=A^2(t) \left(\frac{u^2}{2}\right)_x,\,\,\, q_x =v_x.
\end{aligned}
\]
Using these expressions in, we get
\begin{equation*}
\begin{aligned}
&p_\tau +\left(\frac{p^2}{2} \right)_x=A^2(t) \left(\alpha(t) u+u_t +\left(\frac{u^2}{2}\right)_x \right),\\
 &q_\tau +(p q)_x=A(t)(v_t +(v u)_x)
\end{aligned}
\end{equation*}
Since $A(t) \neq 0$, the result follows.

{ \flushleft
 {\bf Remark:}}  In general 
even for $C^\infty$ initial data with compact support, we have that $u$ is in the space of bounded variation and $v$ is a bounded Borel measure. So it is not obvious to construct the solutions for the system  with the damping term from the equation without the damping term, as the equation is nonlinear. Here we construct solution of initial boundary value problem \eqref{e1.1}  for the case $\alpha=0$.

\subsection{\label{sec:level2}Initial boundary value problem when $\alpha=0$}
In this section, we consider  the system \eqref{e1.1} with $\alpha=0$ i.e. the system for $(p,q)$ with $t=\tau$ namely \eqref{e3.1}
in the quarter plane $x>0,t>0$ with  corresponding initial conditions 
\begin{equation}
\begin{aligned}
p(x,0)=p_0(x),\,\,\,q(x,0)=q_0(x),\,\,\, x>0
\end{aligned}
\label{e3.2}
\end{equation}
and with a weak form of boundary conditions at $x=0$  in
\begin{equation}
\begin{aligned}
p(0,\tau)=p_B(\tau),\,\,\,\int_0^\infty q(y,\tau) dy =q_B(\tau).
\end{aligned}
\label{e3.3}
\end{equation}
Here, we assume $p_0$  is a function of  bounded  variation and $q_0$  an
integrable function which is of  bounded variation in $x >0$.   We also assume   $p_B$ and $q_B$   to be   locally 
Lipschitz continuous functions   in   $\tau>0$.

Indeed
with strong form of the Dirichlet boundary conditions \eqref{e3.3}, there is
neither existence nor uniqueness   of the solution to \eqref{e3.1},  as the speed of propagation
$\lambda =p$ does not have
a definite sign at the boundary $x=0$. We note that the speed is
completely determined by the first equation in \eqref{e3.1}.  We use the  Bardos Leroux and
Nedelec
\cite{Bardos} formulation of the boundary condition for the $u$ component
which for our case is
equivalent to the following condition, see LeFloch \cite{Lefloch}
\begin{equation}
p(0,\tau) \in E(p_B(\tau)),\,\, \text {a.e. } \tau>0.
\label{e3.4}
\end{equation}
For $p_B$ the admissible set $E(p_B)$ is defined by
\begin{equation}
E(p_B) =  \begin{cases} \displaystyle
      {(-\infty,0] ,\,\,\,\text{if} \,\,\,p_B \leq 0}
\\\displaystyle
        {\{p_B\}\cup (-\infty, -p_B],\,\,\,\text{if} \,\,\, p_B>0}.
\end{cases}
\label{e3.5}
\end{equation}
Also $p$ is required to satisfy   the following 
Lax entropy condition \cite{Lax}
\begin{equation*}
p(x-,\tau)\geq p(x+,\tau),\,\, \, a.e. \,\,\, x>0.
\end{equation*}
Also $q$ is required to satisfy the   following  condition on the mass 
\[
\int_0^\infty q(y,\tau)dy =q_B(\tau),
\]
in a weak sense.\\

An explicit solution to a related problem was derived in \cite{Joseph2}. For clarity and completeness we  give the formula for the solution  and give the proof  here. 

There are different  constructions of explicit
representations of the entropy weak solution of the
first component $p$  of \eqref{e3.1} with initial condition \eqref{e3.2}
and the
boundary condition \eqref{e3.4}-\eqref{e3.5}, see  \cite{Joseph3, Joseph5} and  LeFloch \cite{Lefloch}.  The main idea is to use the formula  in \cite{Joseph5} for $p$. This  formula 
involves a
minimisation of functionals on certain class of paths
and generalised characteristics. Once we have $p$,  it is easy to see that the equation for $q$ is a linear equation with discontinuous coefficient.  We get a formula for $q$,  using the fact that the coefficient $p$ is an entropy solution.
First note that  the condition $\int_0^\infty q(y,\tau) dy=q_B(\tau)$ is prescribed only if
the characteristics at $(0,t)$ has positive speed, i.e.  if   $(p(0+,\tau))>0$.
So the weak form of the mass condition for the $q$ component 
is  given by:
\begin{equation}
\text{if $p(0+,\tau)>0$,  then $\int_0^\infty q(y,\tau)dy=q_B(\tau)$.}
\label{e3.6}
\end{equation}
The condition \eqref{e3.6} is required to be satisfied for almost 
every point of the set $\{\tau>0 : p(0+,\tau)>0\}$.
 
 First we introduce a functional defined on 
certain  class of paths,  to describe the solution.
For each fixed $(x,y,\tau), x \geq 0, y \geq 0, \tau>0$,
 the following class of paths $\beta$ in 
the quarter plane
$D=\{ (z,s) : z\geq 0, s \geq 0\}$   is denoted by $C(x,y,\tau)$.   Each path connects  the 
point
$(y,0)$ to $(x,\tau)$ and is of the form $z=\beta(s)$, where $\beta$ is a 
piecewise linear function of maximum three lines. We 
define a functional  on $C(x,y,\tau)$, 
\begin{equation*}
J(\beta) = -\frac{1}{2}\int_{\{s:\beta(s)=0\}} (p_B(s)^{+})^2 ds 
+ \frac{1}{2} \int_{\{s:\beta(s) \neq 0\}} 
\left(\frac{d\beta(s)}{ds}\right)^2 ds.
\end{equation*}

 $\beta_0$  denotes  the straight line path connecting $(y,0)$ and $(x,\tau)$ 
which does not touch the boundary space boundary $x=0$, namely  the space 
$\{(0,t), \tau>0\}$. Let
\begin{equation*}
 A(x,y,\tau)= J(\beta_0) = \frac{(x-y)^2}{2\tau}.
\end{equation*} 
Any $\beta \in C^{*}(x,y,\tau) = C(x,y,\tau)-\{\beta_0\}$
is made up of three pieces lines
connecting $(y,0)$ to $(0,\tau_1)$ in the interior, $(0,\tau_1)$ to 
$(0,\tau_2)$ on the boundary and $(0,\tau_2)$ to $(x,\tau)$ in the interior.
 It can be easily seen that  for such curves,
\begin{equation*}
J(\beta) = J(x,y,\tau,\tau_1,\tau_2) = 
-\int_{\tau_1}^{\tau_2}\frac{(p_B(s)^{+})^2}{2}ds + 
\frac{y^2}{2 \tau_1} + \frac{x^2}{2(\tau-\tau_2)}.
\end{equation*}
It was proved in \cite{Joseph5}, that there exists $\beta^{*} \in 
C^{*}(x,y,\tau)$ 
and corresponding $\tau_1(x,y,\tau),\tau_2(x,y,\tau)$ such that
\begin{equation*}
\begin{aligned}
B(x,y,\tau)&= \min \{J(\beta) :\beta \in C^{*}(x,y,\tau)\}\\
        &= \min \{J(x,y,\tau,\tau_1,\tau_2): \,\, 0\leq \tau_1 < \tau_2 < \tau\}\\
        &=J(x,y,\tau,\tau_1(x,y,\tau),\tau_2(x,y,\tau))
\end{aligned}
\end{equation*} 
is Lipshitz continuous.  For  $P_0$  defined by
\[
 P_0(y)=\int_0^y p_0(s) ds,
 \]
 the functions
\begin{equation*}
\begin{aligned}
Q(x,y,\tau)&= \min\{J(\beta) : \beta \in C(x,y,\tau)\}\\
        & = \min \{A(x,y,\tau),B(x,y,\tau)\}
\end{aligned}
\end{equation*}
and
\begin{equation}
P(x,\tau)= \min \{Q(x,y,\tau) +  P_0(y),\,\,  0\leq y< \infty\}
\label{e3.11}
\end{equation}
 are locally Lipshitz continuous functions in their variables $x, y, \tau.$ Further,  the minimum in \eqref{e3.11}
is attained at some value of $y\geq0$. $y$ depends on $(x,\tau)$, we
call it $y(x,\tau)$. If $A(x,y(x,\tau),\tau)\leq B(x,y(x,\tau),\tau)$
\begin{equation}
\begin{aligned}
P(x,\tau)= \frac{(x-y(x,\tau))^2}{2\tau} +P_0(y(x,\tau))
\end{aligned}
\label{e3.12}
\end{equation}
and if $A(x,y(x,\tau),t)> B(x,y(x,\tau),\tau)$
\begin{equation}
P(x,\tau)=J(x,y(x,\tau),\tau,\tau_1(x,y(x,\tau),\tau),\tau_2(x,y(x,\tau),\tau))+U_0(y(x,\tau))
\label{e3.13}
\end{equation}

Here and hence forth, $y(x,\tau)$ is a minimizer in \eqref{e3.11} and in 
the case of \eqref{e3.12}, \eqref{e3.13},  $\tau_2(x,t)=\tau_2(x,y(x,\tau),\tau)$ and
$\tau_1(x,t)=\tau_1(x,y(x,\tau),\tau)$.
With these notations, we have an explicit formula for the solutions of
\eqref{e3.1} and \eqref{e3.2} with boundary conditions \eqref{e3.4} and
\eqref{e3.6} in the following result.
{ \flushleft
{\bf Theorem 3.2: }}  Assume that $p_0$  is a function of  bounded  variation, $q_0$  an
integrable function which is of  bounded variation in $x >0$,  and  $p_B$, $q_B$   are  locally 
Lipschitz continuous functions   in   $\tau>0$. Then for every $(x,\tau)$,  the   minimum in \eqref{e3.11}
 is achieved by some $y(x,\tau)$ (this $y$ may not be unique) and $P(x,\tau)$ is a 
Lipschitz continuous. 
Further, for almost every $(x,\tau)$ there exists a unique minimum $y(x,\tau)$.  We have  either 
$A(x,y(x,\tau),\tau)< B(x,y(x,\tau),\tau)$ or $B(x,y(x,\tau),\tau) \leq A(x,y(x,\tau),\tau)$. Define
\begin{equation}
p(x,\tau) = \begin{cases}
 \frac{x-y(x,\tau)}{
 \tau},&\text{if } A(x,\tau)< B(x,\tau),\\
 \frac{x}{\tau-\tau_2(x,\tau)},&\text{if } A(x,\tau)> B(x,\tau),
\end{cases}
\label{e3.14}
\end{equation}
and
\begin{equation}
V(x,\tau) = \begin{cases}
- \int_{y(x,\tau)}^\infty q_0(z)dz,&\text{if } A(x,\tau)< B(x,\tau),\\
 - q_B(\tau_2(x,\tau)),
  &\text{if } A(x,\tau)>B(x,\tau). 
\end{cases}
\label{e3.15}
\end{equation}
 Set
\begin{equation}
\begin{gathered}
q(x,\tau)= \partial_{x}(V(x,\tau)).
\end{gathered}
\label{e3.16}
\end{equation}
Then  we have,  the function $(p(x,\tau),q(x,\tau))$ given by 
\eqref{e3.14}-\eqref{e3.16} is a 
weak solution 
of \eqref{e3.1} with initial conditions 
\eqref{e3.2}, boundary conditions \eqref{e3.4}-\eqref{e3.5} and \eqref{e3.6} and the entropy condition $p(x-,\tau) \geq p(x+,\tau)$, a.e. $x>0$.
{ \flushleft
{\bf Proof :}}
First, we recall  from \cite{Joseph3, Joseph5} some properties of
minimizers $y(x,\tau)$  in \eqref{e3.11} and corresponding
$\tau_2(x,\tau)$  and $\tau_1(x,\tau)$ that are required for our
analysis. These minimizers $y(x,\tau)$ may not be unique, for every $(x,\tau)$.
Let
$y^{-}(x,\tau)$ and $y^{+}(x,\tau)$ be the smallest and the largest of these
minimizers  respectively. For each $\tau>0$, they are non-decreasing
functions of
$x$ and hence except for a countable number of points, they
are equal. Corresponding  $\tau_2^{-}(x,\tau)$ and $\tau_2^{+}(x,\tau)$ have the
following properties. They are non-decreasing function of
$x$ for each fixed $\tau$ and hence except for a countable number of
points $x$, they are equal. Also they are  non-decreasing functions of
$\tau$ for each fixed $x$ and except for a countable number of
points $\tau$, they are equal. If
$A(x,y(x,\tau),\tau)<B(x,y(x,\tau),\tau)$ for some $x=x_0$, then this continues
to be so for all $x<x_0$ and $\tau_1(x,\tau)=\tau_1(x_0,\tau)$ when $x\leq x_0$. 
Further, if $A(x,y(x,\tau),\tau)>B(x,y(x,\tau),\tau)$ for some
$x=x_0$, then this continues to be so for all $x>x_0$.  

It was proved in \cite{Joseph5} that,
$p(x,\tau)=Q_1(x,y(x,\tau))=\partial_{x}P(x,\tau)$ is the weak solution of
\begin{equation*}
p_t+\left( \frac{p^2}{2} \right)_x=0
\end{equation*}
with  initial condition $p(x,0)=p_0(x)$, weak form of
boundary condition \eqref{e3.4} and  satisfying  the entropy condition 
$p(x-,\tau)\geq p(x+,\tau)$ for $x>0$,  where
$Q_1(x,y,\tau)= \partial_x Q(x,y,\tau).$ 

To show that $q$ satisfies the second equation, we follow the work of  LeFloch
\cite{LeFloch} and use the
non-conservative product of Volpert \cite{Volpert} in sense of measures.
First, we note that $p$ and $V$ are functions of bounded  
variation and $q=V_x$ is a Borel measure. To justify the product
$p v=p V_x$, we write
\[
[0,\infty)\times[0,\infty) =S_c \cup S_j \cup S_n
\]
where $S_c$ and $S_j$ are points of approximate continuity of $p$ 
and
points of approximate jump of $p$  respectively. $S_n$ is a set of one 
dimensional
Hausdorff-measure zero. At any point $(x,\tau) \in S_j$, $p(x-0,\tau)$ and
$p(x+0,\tau)$ denote the left and right values of $p(x,\tau)$  respectively. For any
continuous function $g :R^1 \to R^1$, the Volpert product
$g(p)V_x$ is defined as a Borel measure in the following manner.
Consider the averaged superposition of $g(p)$, see the work of  Volpert \cite{Volpert},
\begin{equation}
\overline{g(p)}(x,\tau) = \begin{cases}
 g(p(x,\tau)),&\text{if } (x,\tau) \in S_c,\\
  \int_0^1 g(1-\alpha)(p(x-,\tau)+\alpha p(x+,\tau))d\alpha,
  &\text{if }(x,\tau) \in S_j
\end{cases}
\label{e3.18}
\end{equation}
and the associated measure
\begin{equation}
[g(p)V_x](A)=\int_{A}\overline{g(p)}V_x
\label{e3.19}
\end{equation}
where $A$ is a Borel measurable subset of $S_c$ and
\begin{equation}
[g(p)R_x](\{(x,\tau)\})=\overline{g(p)}(x,\tau)(R(x+,\tau)-R(x-,\tau))
\label{e3.20}
\end{equation}
provided $(x,\tau) \in S_j$.

Now, we take the continuity equation.
First, we show that $V$ satisfies
\begin{equation}
\mu = V_\tau+ p V_x=0
\label{e3.21}
\end{equation}
in the sense of measures.
Let $(x,\tau)\in S_c$ and $u=\frac{x-y(x,\tau)}{\tau}$, since $p$
satisfies   \eqref{e3.1}  we have
\[
\begin{aligned}
 \left\{ -\frac{(x-y(x,\tau))}{\tau^2}-\frac{\partial_{\tau}y(x,\tau)}{\tau}+
\left( \frac{x-y(x,\tau)}{\tau} \right) \frac{(1-\partial_{x}y(x,\tau))}{\tau} \right\}=0.
\end{aligned}
\]
 From this, it follows that
\begin{equation}
\partial_{\tau}y(x,\tau) + p \partial_{x}y(x,\tau)=0.
\label{e3.22}
\end{equation}

Now if $A(x,\tau)\leq B(x,\tau)$,  we have 
\[
\partial_{\tau}V(x,\tau) + p\partial_{x}V(x,t)=v_0 (y(x,t)
\{\partial_{\tau}y(x,\tau) + p \partial_{x}y(x,\tau)\}
\]
and from \eqref{e3.22}, we get
\begin{equation}
\partial_{t}V(x,\tau) + p \partial_{x}V(x,\tau)=0. \label{e3.23}
\end{equation}
Now consider the points for which $A(x,\tau)>B(x,\tau)$. If
$(x,\tau)\in S_c$ then $p(x,\tau)=(\frac{x}{\tau-\tau_2(x,\tau)})$, 
we can show that
\[
\partial_\tau(\tau_2(x,\tau)) + p(x,\tau)\partial_x(\tau_2(x,\tau)) =0
\]

and hence
\begin{equation}
\partial_{\tau}V(x,\tau) + p \partial_{x}V(x,\tau)=0, \label{e3.24}
\end{equation}
So for any Borel subset $A$ of $S_c$, it follows from \eqref{e3.23} 
and \eqref{e3.24} that
\begin{equation}
\mu(A)=0.
\label{e3.25}
\end{equation}
Now we consider a point $(s(\tau),\tau) \in S_j$, then
\[
\frac{ds(\tau)}{d\tau} =
\frac{u(s(\tau)+,\tau)+u(s(\tau)-,\tau)}{2}
\]
is the speed of propagation of the discontinuity at this point and
\begin{equation}
\begin{aligned}
&\mu\{(s(\tau),\tau)\}\\
&=-\frac{ds(\tau)}{d\tau}(V(s(\tau)+,\tau)-V(s(\tau)-,\tau))\\
&+\int_0^1 (p(s(\tau)-,\tau)
+\alpha (p(s(\tau)+,\tau) - p(s(\tau)-,\tau))d\alpha (V(s(\tau)+,\tau)-V(s(\tau)-,\tau))\\
&=\left[-\frac{ds(\tau)}{d\tau} +
\frac{p(s(\tau)+,\tau))+u(s(\tau)-,\tau)}{2} \right]
(V(s(\tau)+,\tau)-V(s(\tau)-,\tau))\\
&=0.
\end{aligned}
\label{e3.26}
\end{equation}
From \eqref{e3.25}, \eqref{e3.26} and \eqref{e3.18}-\eqref{e3.20},  equation \eqref{e3.21} follows. Differentiating \eqref{e3.21} w.r.t. $x$ and
using $q = V_x$, we get $q_\tau +(p q)_x =0$ in the sense of 
distribution.\\

To show that the solution satisfies the initial conditions, first
we observe that given $\epsilon>0$ there exists $\delta>0$ such that for
all $x\geq\epsilon$, $\tau\leq \delta$,
\[
p(x,\tau)=\left( \frac{x-y(x,\tau)}{\tau} \right),\,\,\,V(x,\tau)=-\int_{y(x,\tau)}^\infty q_0(y)dy
\]
 where $y(x,\tau)$  minimizes
$\min_{y\geq 0}[P_0(y)+(\frac{(x-y)^2}{\tau})]$. 
Then
Lax's argument \cite{Lax} gives $\lim_{\tau\to 0}p(x,\tau)=p_0(x)$
for a.e.
$x\geq \epsilon$. Since $\epsilon>0$ is arbitrary,
\begin{equation*}
\lim_{\tau \to 0}p(x,\tau)=p_0(x),\,\,\, \text{ for a.e. }\,\,x.
\end{equation*}
Since $y(x,\tau)-x = -\tau u(x,\tau)$, it follows that
$y(x,\tau)\to x$ as $\tau \to 0$ for a.e $x$. Since the map $y \rightarrow 
-\int_y^\infty q_0(z) dz$ is continuous, we get
\[
\lim_{\tau \to 0}V(x,\tau)=\lim_{\tau\to
0}V(x,\tau))=-\int_x^\infty q_0(z) dz \,\,\, \text{ for a.e. }\,\,x.
\]

Now we show the solution satisfies the boundary condition \eqref{e3.4}
and \eqref{e3.6}.  The fact that the $p$ component
satisfies the boundary condition \eqref{e3.4} is proved in \cite{Joseph5}.
 We need  to show that $q$ satisfies the boundary condition \eqref{e3.6}. Note that
if $p(0+,\tau_0)>0$, then by the weak form of the boundary condition for $p$, $p(0,\tau_0)=p_B(\tau_0)>0$. Then we have that $p_B(\tau)>0$ in a neighbourhood of $\tau_0$. By the method of characteristics, we can construct the  solution in a small neighbourhood of $(0,\tau_0)$ in $x>0, t>0$ where $p(x,\tau)$ is Lipschitz continuous. Also  $p(x,\tau)>0$ for $0<x\leq \epsilon$
for some sufficiently small $\epsilon$ and $p$ is given by
\[
p(x,\tau)=\frac{x}{\tau-\tau_2(x,\tau)}.
\]
We have $\tau-\tau_2(x,\tau) = x/(p(x,\tau))$ and it follows that
$\lim_{x\to 0}\tau_2(x,\tau)=\tau$. So we have 
\[
\int_0^\infty q(x,\tau)dx=-V(0+,\tau)=\lim_{x\to 0}q_B(\tau_2(x,\tau))=q_B(\tau). 
\]
The proof of the theorem is complete.

\subsection{\label{sec:level3}Solution to the initial boundary value problem when $\alpha(t) \neq 0$}

Under the transformation \eqref{3.3}, equation \eqref{e1.1} transforms to equation \eqref{e3.2}, the initial conditions \eqref{e1.2} become \eqref{e3.3} with $(p_0,q_0) =(u_0,v_0)$. The boundary conditions \eqref{e1.3} become $p(0,\tau)=p_B(\tau)=A(t(\tau))u_B(t(\tau)$,  $\int_0^\infty q(y,\tau)dy=q_B(\tau)=v_B(t(\tau))$. With these initial and boundary values
in previous theorem we formally get a formula for the solution in the following theorem. 

{ \flushleft
{\bf Theorem 3.3:}}  Assume  that $u_0$ is a function of bounded variation, $v_0$ is integrable and a function of bounded variation.  Also assume that  $u_B(t)$ and $v_B(t)$ are locally Lipschitz continuous functions on $[0,\infty)$.
Let $(p(x,\tau),q(x,\tau))$ be solution given by  Theorem 3.2 in \eqref{e3.14}-\eqref{e3.16} with initial condition \eqref{e3.2} and boundary conditions \eqref{e3.4} and \eqref{e3.6}, then
$(u(x,t),v(x,t))$ defined by
 \begin{equation}
\begin{aligned}
 u(x, t) =   e^{-\int^t \alpha(s) ds} p(x,\tau(t))  , \,\,\, v(x, t)= q(x,\tau(t))
\end{aligned}
\label{e3.29}
\end{equation}
is a weak solution of  \eqref{e1.1}  satisfying  initial condition \eqref{e1.2}  and  a weak form of boundary condition \eqref{e1.3} , namely $u(0+,t)$ satisfies \eqref{e1.7}- \eqref{e1.8} and if $u(0+,t)>0$,
then, $\int_0^\infty v(y,t) dy =v_B(t)$, a.e. $x>0$.
{ \flushleft
{\bf Proof :}} We need to verify that $(u,v)$ satisfies  equation \eqref{e1.1}. Let us denote $\bar{V}(x,t)=V(x,\tau(t))$, then $v(x,t)=q(x,\tau(t))=\bar{V}(x,t)$. 
Since $\tau(t)=\int_0^t e^{-\int_0^s \a(s')ds'} ds$, by chain rule we have
\[
\begin{aligned}
&u_t=-\a(t) u(x,t) +e^{-2 \int_0^t \a(s')ds'} p_\tau, \,\,\,u_x=e^{-2 \int_0^t \a(s')ds'} p_x\\
&\bar{V}(x,t)_t= V_\tau(x,\tau) \frac{d\tau}{dt} =V_\tau(x,\tau(t)) e^{-\int_0^t \a(s')ds'},\,\, \bar{V}(x,t)_x= V_x
\end{aligned}
\]
These relations leads to 
\[
\begin{aligned}
&u_t +u u_x +\a(t) u =e^{-2\int_0^t \a(s')ds'} (p_\tau +p p_x)=0\\
&\bar{V}_t +u \bar{V}_x=e^{-\int_0^t \a(s')ds'} (V_\tau +p V_x)=0
\end{aligned}
\]
where we used the fact that $(p,q)(x,\tau)$ satisfies the equation for $\a=0$.
The proof that $(u,v)$ given by \eqref{e3.29} satisfies the initial condition \eqref{e1.2} and boundary conditions \eqref{e1.3} in the weak sense, follows exactly along the same arguments as in the case when $\alpha=0$.\\

\subsection{The Initial Value Problem} 
In this section we consider the initial value problem for \eqref{e1.1}  in $x \in R^1, t>0$ with initial conditions 
\begin{equation}
 \begin{aligned}
u(x, 0)= u_0(x), \,\,\, v(x, 0)= v_0(x),\,\,\, x \in R,
\end{aligned}
\label{e3.30}
\end{equation}

 We have the following theorem.  The explicit formula derived here agrees with the one obtained in \cite{Divya}, by the vanishing viscosity approximation \eqref{e2.1} for the special case when $\alpha$ is a constant.
{ \flushleft
{\bf Theorem 3.4:}} Suppose $(u_0,v_0)$ is a pair of functions of bounded variation and $v_0$ is integrable in $ x \in R$. Let $y(x,\tau)$ be a minimiser for 
\[
\min_{-\infty <y<\infty} \left\{ \frac{(x-\tau)^2}{2 \tau} + \int_0^y u_0(z) dz \right\}
\]
which exists for all $(x,t), x \in R, t>0$ and is unique  a.e.  Then the distribution given by 
\begin{equation*} 
(u(x,t),v(x,t))= \left( \frac{1}{A(t)} \left( \frac{x-y(x,\tau(t))}{\tau(t)} \right),- \partial_x \left(\int_{y(x,\tau(t))}^\infty v_0(z)dz \right) \right),
\end{equation*}
is a solution to the initial value problem for \eqref{e1.1}  with initial conditions \eqref{e3.30}. \\

The proof is exactly same as in Theorem 3.3, with the only difference  that the solution to the initial value problem for  $\alpha=0$  in \cite{Joseph4,Joseph7} must be used. \\

Similarly, we can get an explicit formula for solutions of \eqref{e1.1} in $x>0,t>0$ with initial conditions  \eqref{e1.2}  and  a weak form of the Dirichlet boundary conditions  \eqref{e1.4}. 
instead of the mass condition on $v$. Here again after transforming the problem to the case when $\alpha=0$  we may use the formula in \cite{Joseph6}  to  get an explicit formula for the solution. The details are omitted.

\section{Conclusion}

 The main difficulty in studying  the system  for the  Eulerian droplet model for air particle flow \eqref{e1.1}  is that the solutions do not belong to the classical $L^\infty$ or the  space of functions of bounded variation but generally are in the space of measures, even if we start with  smooth initial data. This is due to concentration of mass which leads to the formation of $\delta$ waves in the  $v$ component. The same phenomenon happens when $\alpha= 0$  also.  The significance of our work is  that we consider the case $\alpha(t) \neq 0$  and   we also consider the initial  boundary value problem with   general initial and and boundary data.  In this work  we show the existence of the solution to the Cauchy problem for the parabolic approximation of \eqref{e2.1} with general  initial condition \eqref{e2.2}  and boundary conditions \eqref{e2.3}  by adding a viscous term with a small coefficient.  Using this result, we prove the existence of a  weak asymptotic solution  for \eqref{e1.1} having general  initial and boundary data. Next we  use a transformation to reduce the system \eqref {e1.1} to the case when  $\alpha=0$.  Then we derive an explicit formula for solution  with the velocity  $u$ in the space of functions of bounded variation and volume fraction $v$ in the space of bounded Borel measures. Here, we need a weak formulation of the boundary condition.  This construction is done using the Hopf-Lax formula  for the initial boundary value problem for the Burgers equation and the non-conservative product of Volpert  \cite{Volpert}.

\begin{acknowledgments}
{ \flushleft
 
This work is done while the author was  at IISER Pune as postdoctoral fellow.  The author would like to thank the referee for the constructive comments and suggestions that has improved this paper. 
}
\end{acknowledgments}

\section*{AUTHOR DECLARATIONS}
{ \flushleft
{ \bf Conflict of interest disclosure:} None. \\
\vspace{0.25 cm}
{ \bf Ethics approval statement:}  Not applicable. \\
\vspace{0.25 cm}
{ \bf Author contributions.} \\
{\bf Kayyunnapara Divya Joseph:} Formal analysis (complete work); Investigation (complete work). \\
\vspace{0.25 cm}
{ \bf Funding statement:}  No funding involved.  

\section*{Data Availability Statement}
 Not applicable. 
}
\section*{ORCID}
Kayyunnapara Divya Joseph https://orcid.org/0000-0002-4126-7882

\bibliography{aipsamp}% Produces the bibliography via BibTeX.

\end{document}